\documentclass{amsart}

\newtheorem{theorem}{Theorem}[section]
\newtheorem{lemma}[theorem]{Lemma}

\theoremstyle{definition}
\newtheorem{definition}[theorem]{Definition}
\newtheorem{example}[theorem]{Example}

\theoremstyle{remark}
\newtheorem{question}[theorem]{Question}

\theoremstyle{remark}
\newtheorem{remark}[theorem]{Remark}

\numberwithin{equation}{section}

\usepackage[latin1]{inputenc}
\usepackage{amsfonts}
\usepackage{amssymb}
\usepackage{amsmath}
\renewcommand{\phi}{\varphi}
\newcommand{\C}{{\mathbb{C}}}
\newcommand{\R}{{\mathbb{R}}}
\newcommand{\Q}{{\mathbb{Q}}}
\newcommand{\Z}{{\mathbb{Z}}}
\newcommand{\N}{{\mathbb{N}}}

\renewcommand{\epsilon}{\varepsilon}
\renewcommand{\theta}{\vartheta}

\newcommand{\w}{\wedge}

\DeclareMathOperator{\rk}{rk}

\DeclareMathOperator{\lcm}{lcm}
\DeclareMathOperator{\diag}{diag}

\begin{document}

\title[Contact homology of Brieskorn manifolds]{Contact homology of Brieskorn manifolds}

%    Information for first author
\author{Otto van Koert}
\address{Otto van Koert\\
Universit\'e Libre de Bruxelles\\
D\'epartement de Math\'ematiques - CP 218\\
Boulevard du Triomphe\\
B-1050 Bruxelles - Belgique}
\email{ovkoert@ulb.ac.be}


\begin{abstract}
  We give an algorithm for computing the contact homology of some Brieskorn
  manifolds. Brieskorn manifolds can be regarded as circle-bundles over
  orbifolds and our algorithm expresses cylindrical contact homology of the
  Brieskorn manifold in terms of
  the homology of the underlying orbifold. As an
  application, we construct infinitely many contact structures on the class
  of simply connected contact manifolds that admit nice contact forms
  (i.e.~no Reeb orbits of degree -1, 0 or 1) and have index positivity with
  trivial first Chern class.
\end{abstract}

\maketitle

\section{Introduction}
For a long time it has been known that contact manifolds may carry many
non-isomorphic contact structures. A first way to distinguish these
structures from each other is by considering their Chern class or their
formal homotopy class. In dimension 3 we can, in addition, sometimes
distinguish contact structures by showing tightness or overtwistedness as
was shown first by Bennequin. At present, the latter two notions do not
have generalizations to higher dimensions. 

Eliashberg was the first to show that the classical invariants (i.e.~the
Chern class and formal homotopy class) are not always enough for a full
classification of contact structures on a given manifold in dimension
larger than 3. He showed that spheres in dimensions $4n+1$ admit a
non-standard, yet homotopically trivial contact structure in
\cite{Eliashberg_filling}. Giroux gave an interesting example in dimension
3 in \cite{Giroux_torus}. He exhibited infinitely many tight contact
structures on $T^3$ in the same homotopy class of plane fields. Since then more
general techniques to distinguish contact structures were introduced.
Eliashberg and Hofer's contact homology is such a technique (for a survey,
see \cite{Eliashberg_survey}, and for a more recent description,
\cite{Eliashberg}). It works, roughly speaking, as follows. The Reeb orbits
on a given contact manifold can be given a grading. The chain complexes are
freely generated by the closed Reeb orbits. Then a differential can
be defined by counting certain holomorphic curves asymptotic to Reeb orbits
in the symplectization of the contact manifold. It can be shown that the
homology of this differential is independent of the contact form under
suitable conditions and it is hence an invariant of the contact structure.
Since its introduction, several examples of non-isomorphic contact
structures in the same homotopy class have been found.

One of these examples is a family of Brieskorn spheres in dimensions
$4n+1$. These were found by Ustilovsky in \cite{Ustilovsky}. In his thesis he
used a few other tools that can easily be applied to construct more contact
structures on the same manifold. This is done by connect summing. It is
well known that the connected sum of two contact
manifolds carries again a contact structure, see \cite{Meckert} and \cite{Weinstein}. In his
thesis Ustilovsky showed that we can roughly see what happens to the
resulting contact homology, provided we are given sufficient information on
the given two contact manifolds.

In this paper, we use these ideas to construct more contact structures on a
certain class of contact manifolds. The starting point is the computation
of the contact homology of Brieskorn manifolds. This is done by calculating
the Morse-Bott contact homology, an extension of contact homology due to
Bourgeois \cite{Bourgeois_thesis}. It allows a larger class of contact
forms to be admissible than the generic contact forms required in the
original formulation of contact homology (\cite{Eliashberg}). Using
Bourgeois's techniques we give an algorithm that computes the contact
homology of a Brieskorn manifold. 

Notice that a Brieskorn manifold can be regarded as an $S^1$-bundle over an
orbifold. It turns out that the resulting contact homology of the Brieskorn
manifold with its natural contact structure can be expressed in terms of
the homology of that orbifold and its singular subspaces with some degree
shifts that we explicitly compute. By combining this algorithm with
Ustilovsky's ideas we produce a large class of contact manifolds with
infinitely many non-isomorphic contact structures.

\subsubsection*{Acknowledgements}
This paper is a part of my thesis, which I wrote under supervision of H.~Geiges at the university of Cologne. I am indebted to him for his support and patience. I am
very grateful to F.~Bourgeois for his explanations and for many valuable
comments and suggestions. I would also like to thank R.~Lipshitz for his comments. Currently I am supported by the F.N.R.S., Belgium.

\section{Preliminaries}\label{prelims}
The \textbf{Brieskorn manifolds} $\Sigma (a_0,\ldots,a_n)\subset\C^{n+1}$ (with
$a_0,\ldots,a_n\in\N$) are defined as the intersection of the sphere
$S^{2n+1}=\{(z_0,\ldots,z_n)\in \C^{n+1}~ |~ |z_0|^2+\ldots+|z_n|^2=1\}$ with the
zero set of the polynomial
$f(z_0,\ldots,z_n)=z_0^{a_0}+\ldots+z_n^{a_n}$. The Brieskorn manifolds
carry a contact structure which was shown by Lutz and Meckert
\cite{LutzMeckert}. We will take the following contact form on $\Sigma (a_0,\ldots,a_n)$:
$$
\alpha=\frac{i}{8}\sum_{j=0}^n a_j(z_j \mathrm d \bar z_j-\bar z_j \mathrm d z_j).
$$
We will denote the contact structure by $\xi=\ker \alpha$.
The Reeb vector field of this contact form has a particularly simple shape
$$
R=4i (z_0/a_0,\ldots,z_n/a_n),
$$
where we regard $T \Sigma (a_0,\ldots,a_n)$ as a subset of $T\C^{n+1}$. The
Reeb flow then looks like
\begin{equation}
\label{eq:Reebflow}
\phi_t(z)=(e^{4it/a_0}z_0,\ldots,e^{4it/a_n}z_n),
\end{equation}
which means that all Reeb orbits are closed. In particular, all Reeb orbits are degenerate. We can, however, still compute the contact
homology using the Morse-Bott approach due to Bourgeois
\cite{Bourgeois_thesis}, without an explicit perturbation of the contact
form. In the rest of this section we will establish some notation and set up
some ingredients of the Morse-Bott approach.

The first step will be the computation of the Maslov indices of the Reeb
orbits. They will give the grading of the Morse-Bott complex. To compute the Maslov indices we will be using
the same method as Ustilovsky did in \cite{Ustilovsky}: we extend the Reeb field
so that it gives rise to a symplectomorphism of $\C^{n+1}$. We may then do the
computations in $\C^{n+1}$. In our discussion, the symplectic form on $\C^{n+1}$ is given by
$$
\omega=\mathrm d\alpha=\frac{i}{4}\sum_{j=0}^n a_j \mathrm d z_j \wedge \mathrm d \bar z_j. 
$$ 
We start by showing that we can construct a basis of the symplectic
complement $\xi^{\omega}$ of $\xi$ in $\C^{n+1}$. Notice that
$\xi^{\omega}=\mathrm{span}(X_1,Y_1,X_2,Y_2)$, where
\begin{eqnarray*}
 X_1=(\bar z_0^{a_0-1},\ldots,\bar z_n^{a_n-1}), &~& Y_1=iX_1, \\
 X_2=-2i(\frac{z_0}{a_0},\ldots,\frac{z_n}{a_n}), &~& Y_2=(z_0,\ldots,z_n).
\end{eqnarray*}
We are going to use a ``Gram-Schmidt'' process to turn this into a
symplectic standard basis. This new basis is given by
\begin{eqnarray*}
\tilde X_1 =\frac{X_1}{\sqrt{\omega(X_1,Y_1)}}, &~&\tilde Y_1
=\frac{Y_1}{\sqrt{\omega(X_1,Y_1)}}  \\
\tilde X_2=X_2, &~&\tilde Y_2=Y_2-
\frac{\omega(X_1,Y_2)Y_1-\omega(Y_1,Y_2)X_1}{\omega(X_1,Y_1)}=Y_2-\frac{\sum
  a_i z_i^{a_i}}{2\omega(X_1,Y_1)}X_1,
\end{eqnarray*}
where we have used that $\omega(X_1,Y_1)=\frac{1}{2} \sum_j a_j
|z_j|^{2(a_j-1)}>0$. This new basis is a standard symplectic basis of
$\xi^{\omega}$. Note that $\xi \oplus \xi^{\omega}=\C^{n+1}$.

This computation gives us a bonus: $c_1(\xi)=0$. Indeed, both
$\xi^{\omega}$ and $\C^{n+1}$ are symplectically trivial and have therefore
trivial total Chern class. Hence the Chern class of $\xi$ is trivial as well.

\subsection{Maslov indices}\label{maslov_section}
The Reeb flow $\phi$ which we introduced in the previous section can obviously be extended
to a flow on $\C^{n+1}$, which we will also denote by $\phi$. We find the
following for the action of the extended Reeb flow on $\xi^\omega$:
\begin{eqnarray}
\label{eq:maslov_complement}
T\phi_t(\tilde X_1(x))=e^{4it} \tilde X_1( \phi_t(x)), 
&~& T\phi_t(\tilde Y_1(x))=e^{4it} \tilde Y_1( \phi_t(x))
\end{eqnarray}
\begin{eqnarray*}
T\phi_t(\tilde X_2(x))=\tilde X_2( \phi_t(x)), 
&~& T\phi_t(\tilde Y_2(x))=\tilde Y_2( \phi_t(x)).
\end{eqnarray*}
The action of $T\phi_t$ on $\C^{n+1}$ is given by the differential of
$\phi_t$. In terms of the standard basis of $\C^{n+1}$, it is given by the
diagonal matrix
$$
\diag(e^{4it/a_0},\ldots,e^{4it/a_n}).
$$ 
We can now use the additivity of the Maslov index under direct sums of symplectic paths
to get the index for $\xi$ (for general properties of the Maslov index see
\cite{RobbinSalamon} and \cite{Ustilovsky_thesis}). Let $\gamma=\{ \phi_t(p)~|~t\in[0,T] \}$ be a
closed Reeb orbit of period T through $p\in \Sigma (a_0,\ldots,a_n)$. Let
$\Phi_{\C^{n+1}}(t)=T\phi_t(p)$ for $t\in [0,T]$, the path of symplectic
matrices associated to the Reeb flow extended to $\C^{n+1}$. Write
$\Phi_{\xi^{\omega}}(t)=T\phi_t(p)|_{\xi^{\omega}}$ for $t\in [0,T]$, the
path of symplectic matrices induced by the Reeb flow on the symplectic
complement of $\xi$. Then
\begin{equation}
\label{eq:maslov_reeb}
\mu(\gamma)=\mu( \Phi_{\C^{n+1}} )-\mu( \Phi_{\xi^{\omega}} ).
\end{equation}
Note also that the right-hand side of the above equation can be easily
computed using additivity of the Maslov index under direct sums and the
following formula (see \cite{Ustilovsky_thesis}):
\begin{equation}
\label{eq:maslov_unitair}
\mu( e^{it}|_{t \in [0,T]} )=
\left\{
\begin{array}{cc}
\frac{T}{\pi} & \text{if } T \in 2\pi \Z \\
2 \left\lfloor \frac{T}{2\pi} \right\rfloor +1 & \text{otherwise.} 
\end{array}
\right.
\end{equation}

\subsection{Homology of Brieskorn manifolds}\label{Brieskorn_homology}
In \cite{Randell} Randell proves an algorithm that computes the homology of
generalized Brieskorn manifolds. Since we are here only interested in the
homology of Brieskorn manifolds, we list the algorithm in the Brieskorn
case and refer the reader to \cite{Randell} for further details. Let
$M=\Sigma(a_0,\ldots,a_n)$. For convenience, we introduce the following
notation. Let $I$ denote the set $\{ 0,\ldots,n \}$. A subset of $I$ with $s$
elements will by denoted by $I_s$. If $I_s=\{ i_1,\ldots,i_s \}$, then let $K(I_s)$ denote the Brieskorn
manifold $\Sigma(a_{i_1},\ldots,a_{i_s})$ of dimension $2s-3$. Note that $M=K(I)$
contains all manifolds $K(I_s)$ for all $1<s\leq n+1$ in a natural way by
restricting to suitable coordinate hyperplanes. We define following Randell
$$
\kappa(K(I_s))=\sum_{I_t\subset I_s} (-1)^{s-t}
\frac{\prod_{i\in I_t}a_i}{\underset{j\in I_t}{\lcm} a_j}.
$$
Then we have for the free part of the homology 
$$
\rk \tilde H_{n-1}(M,\Z)=\kappa (K(I)).
$$
For the torsion part, a few additional definitions are required. We set 
$$
k(K(I_s))=\left\{ 
\begin{array}{ll} 
\kappa(K(I_s)) & \text{if } n+1-s \text{ is odd} \\
0 & \text{otherwise.}
\end{array}
\right.
$$
Let $C(\emptyset)=\gcd_{i\in I} a_i$ and set
$$
C(I_s)=\frac{\underset{i \in(I-I_s)}{\gcd} a_i}{\prod_{I_t \subsetneqq I_s} C(I_t) }.
$$
Now set $d_j=\prod_{k(K(I_s))\geq j}C(I_s)$ and $r=\mathrm{max}\{ k(K(I_s)) | I_s
\subset I \}$. Then we have
$$
\mathrm{Tor}H_{n-1}(K(I),\Z)=Z_{d_1}\oplus \ldots \oplus Z_{d_r}.
$$

There is another interesting result from Randell's paper that we will
use. Set $d=\lcm (a_0,\ldots,a_n)$ and $q_i=d/a_i$. The Brieskorn manifolds admit the following $S^1$-action:
$$
t(z_0,\ldots,z_n)=(t^{q_0}z_0,\ldots,t^{q_n}z_n)
$$
for $t\in S^1$. Then define $M^*:=M/S^1$, which, in general, will be an
orbifold. \label{S1_remark} It is important to note here that this
$S^1$-action coincides with the flow of the Reeb field. Again from \cite{Randell} we can compute the rational homology of
$M^*$. The result is
\begin{equation}
\label{eq:rathom_orbit}
H_q(M^*,\Q) \cong \left\{ \begin{array}{l} \Q, ~q \textrm{ even } 0
    \leq q \leq \dim M^* \\ 0 \text{, otherwise} \end{array} \right\} \oplus
\left\{ \begin{array}{l} \Q^{\kappa},~q=\frac{1}{2} \dim M^* \\ 0 \text{, otherwise} \end{array}\right\},
\end{equation}
where $\kappa=\kappa(K(I))$.

\subsection{Morse-Bott contact homology}\label{MorseBott_homology} Let $M$
denote a manifold with contact structure $\xi$ which is given by as the
kernel of the 1-form $\alpha$. Define the action functional
\begin{eqnarray*}
\mathcal A: C^{\infty}(S^1,M) & \to & \R \\
\gamma &\mapsto & \int_{\gamma} \alpha.
\end{eqnarray*}
The critical points of this action functional are the periodic orbits of
the Reeb flow $\phi_t$. Cylindrical contact homology can then be regarded
as a kind of Morse theory for this action functional, where one needs
appropriate regularity conditions. For instance one would need a generic
contact form for which all Reeb orbits are non-degenerate. In particular
this means that symmetric contact forms like the one for the Brieskorn
manifolds are excluded. This makes computations in general
quite hard. The Morse-Bott approach to contact homology due to Bourgeois is often
much better suited to symmetric contact forms. We will now summarize
some of Bourgeois's results from \cite{Bourgeois_thesis}. As we need some
of his lemmas later on, we will list those along with one of the main
theorems of \cite{Bourgeois_thesis}. We will restrict in some of his
theorems to contact forms for which all Reeb orbits are closed. This is not
a requirement in Bourgeois's work and only takes care of a small technical detail.

Let $\sigma(\alpha)$ denote the action spectrum of $\alpha$,
i.e.~the critical values of the action functional $\mathcal A$.

\begin{definition}\label{MorseBott_definition}
A contact form $\alpha$ is said to be of Morse-Bott type if the action spectrum
$\sigma(\alpha)$ is discrete and if, for every $T'\in\sigma(\alpha)$,
$M_{T'}=\{ p\in M| \phi_{T'}(p)=p\}$ is a closed, smooth submanifold of $M$, such
that the rank $d\alpha|_{M_{T'}}$ is locally constant and $T_pM_{T'}=\mathrm{ker}(T\phi_{T'}-id)_p$.
\end{definition}

The Reeb flow induces an $S^1$-action on $M_{T'}$. Using this action we define
the orbit space $S_{T'}:=M_{T'}/S^1$. Note that these orbit spaces are
orbifolds in general.

The chains of the Morse-Bott chain complex will correspond to the critical
points of suitable Morse functions on the orbit spaces. Bourgeois
constructs these Morse functions by induction (note that he needs Morse functions on
orbifolds for which he introduces a suitable notion in
\cite{Bourgeois_thesis}). We will now describe his construction.

For the smallest $T\in \sigma(\alpha)$, the orbit space $S_T$ is a
smooth manifold. Take any Morse function $f_T$ on it. For larger $T\in
\sigma(\alpha)$, $S_T$ is an orbifold where the singularities are the orbit
spaces $S_{T'}$ with $T'$ dividing $T$. The previously defined Morse functions
$f_{T'}$ on the orbit spaces $S_{T'}$ are extended to a function $f_T$ on
$S_T$ by requiring that the Hessian of $f_T$ restricted to the normal
bundle of $S_{T'}$ is positive definite. These Morse functions are then lifted to $M_T$ and extended to a function
$\bar f_T$ on $M$ such that they have support only in a tubular
neighborhood of $M_T$. 

For $T\in\sigma(\alpha)$ Bourgeois considers the following family of
contact forms $\alpha_\lambda =(1+\lambda \bar f_T)\alpha$. We have
\begin{lemma}[Bourgeois] \label{Bourgeois_perturb} For all $T$, we can choose $\lambda >0$ small
  enough such that the periodic orbits of $R_{\alpha_\lambda}$ in $M$ of
  action $T'\leq T$ are non-degenerate and correspond to the critical
  points of $f_{T'}$. 
\end{lemma}
Let $p\in S_{T'}$ be a critical point of $f_{T'}$ and denote its corresponding
closed Reeb orbit (and its multiple covers) by $\gamma^p_{kT'}$ for $k=1,2,\ldots$.
As Bourgeois's construction is explicit enough, he is able to compute the
Conley-Zehnder indices of these Reeb orbits for small $\lambda$ such that
Lemma \ref{Bourgeois_perturb} applies:
\begin{equation}
\label{eq:CZindex}
\mu_{CZ}(\gamma^p_{kT'})=\mu(S_{kT'})-\frac{1}{2}
\dim S_{kT'}+\mathrm{index}_p(f_{kT'}).
\end{equation}
This determines the degree of the Reeb orbits associated with the perturbed
contact form with small period. The following notion is helpful in dealing
with orbits of larger period.
\begin{definition}
The orbit spaces $S_T$ are said to have {\bf index positivity} if there exist
constants $c>0$ and $c'$ such that $\mu(S_T)>cT+c'$ for all $T\in\sigma(\alpha)$.
\end{definition}
Similarly, we define {\bf index negativity} of the orbit spaces $S_T$ if
there are constants $c<0$ and $c'$ such that $\mu(S_T)<cT+c'$.
In order to control the behavior of orbits 
with larger period, Bourgeois has the following result.
\begin{lemma}[Bourgeois] \label{Bourgeois_pos}
Assume that the orbit spaces $S_T$ have index positivity, that $c_1(\xi)=0$
and that all Reeb orbits are closed. Then there exists a
$\lambda_0>0$ such that, if $0<\lambda<\lambda_0$, all period orbits
$\gamma_\lambda$ of $R_{\alpha_\lambda}$ of action greater than $T$ satisfy
$\mu_{CZ}(\gamma_\lambda)>\frac{c}{2} T$, where $c$ is the positive constant from the index positivity of $S_T$.
\end{lemma}
This lemma makes sure that other closed Reeb orbits that do not correspond
to any critical point of the Morse functions on the orbit spaces have large
Maslov indices. A similar result holds in case of index negativity of the orbit spaces.

The chains of the Morse-Bott complex are the critical points $p$ of the Morse
functions $f_T$ for $T\in \sigma(\alpha)$ with degree given by
\begin{equation}
\label{eq:MorseBott_degree}
\mathrm{deg}(p)=\mu(S_{kT'})-\frac{1}{2}
\dim S_{kT'}+\mathrm{index}_p(f_{kT'}) +n-3.
\end{equation}

\subsubsection{Moduli spaces of generalized holomorphic curves}
Before we come the definition of the moduli spaces, first recall that the
fibered product of $A$ and $B$ over $C$ for maps $f: A \to C$, $g: B\to C$
is given by 
$$
A\times_C B=\{ (a,b)\in A\times B ~| ~f(a)=g(b) \}.
$$ 
Let $\mathcal M(S^+,S^-)$ denote the moduli space of holomorphic curves
with degenerate asymptotics, with orbit spaces $S^+$ and $S^-$ for the
positive and negative puncture, respectively. In other words, elements of $\mathcal M(S^+,S^-)$ correspond to holomorphic cylinders that are asymptotic to a closed Reeb orbit in $S^+$ at the positive end and asymptotic to a closed Reeb orbit in $S^-$ at the negative end of the cylinder. In special cases, such as simple holomorphic curves (which for instance happens if the orbits in $S^+$ and $S^-$ are $k$- and $l$-fold covers of simple Reeb orbits and $k$ and $l$ are relatively prime) one can choose the almost complex structure in such a way that the moduli space is a smooth manifold. In general the moduli space of holomorphic curves is not a smooth manifold though.

In the Morse-Bott picture we shall make use of the moduli space of generalized holomorphic cylinders given by
$$
\mathcal M^{f_T}(S^+,S^-) = \mathcal M(S^+,S^-) \bigcup \left( \mathcal
M(S^+,S')\times_{S'} \left( \R^+\times \mathcal M (S',S^-) \right) \right)
\bigcup \ldots
$$
where the union runs over successive fibered products. Note that the union
is finite, because a holomorphic curve has to have positive energy (the
energy is equal the action of the top Reeb orbit minus the action of the
bottom Reeb orbit) and the action spectrum is discrete. The projection maps
for these fibered products are $ev^-$ and $\phi^{f_T}\circ ev^+$. The maps $ev^-$
and $ev^+$ are the evaluation maps at the negative and positive puncture,
respectively, i.e.~$ev^{\pm}:\mathcal M(S^+,S^-)\to S^{\pm}$ and here $\phi^{f_T}$
is the gradient flow of $f_T$ on the orbit space. So $\phi^{f_T}\circ ev^+:~(t,u)\mapsto \phi^{f_T}(t) \left( ev^+(u)\right)$,
where $t\in \R^+$ and $u\in \mathcal{M}(S',S'')$. 

In other words, generalized holomorphic cylinders are ``cascades'' of holomorphic cylinders and gradient trajectories. For example, the second term of the union consists of a holomorphic cylinder starting asymptotically at a closed Reeb orbit in $S^+$ and ending (not asymptotically) at a closed Reeb orbit in $S'$. After that the cascade follows the gradient flow of the Morse function $f_T$ on the orbit space $S'$ for a finite time. Finally the cascade continues with a holomorphic cylinder ending asymptotically at a closed Reeb orbit in $S^-$. For a discussion of the general case (i.e.~not only cylinders) see \cite{Bourgeois_thesis}, page~38.

\subsubsection{Orientation}
The moduli space of holomorphic curves needs to be oriented. Bourgeois does
this using the coherent orientation scheme, in a similar way as in the original
formulation of contact homology. There are a few obstructions to the
existence of coherent orientations. First of all, there is the notion of a
bad Reeb orbit. Since we are working here with Morse-Bott instead of
generic contact forms, the notion is slightly different from the standard
one described in \cite{Eliashberg}.
\begin{definition}
  A Reeb orbit $\gamma$ is said to be {\bf bad} if it is the $2m$-fold cover of a
  simple Reeb orbit $\gamma'\in S_T$ and if 
$$(\mu(S_{2T})\pm\frac{1}{2}
  \dim S_{2T})-(\mu(S_{T})\pm\frac{1}{2} \dim S_{T})$$ is odd.
  If a Reeb orbit is not bad, it is said to be {\bf good}.
\end{definition}
In Morse-Bott contact homology, there can be another obstruction for
coherent orientations. For this, we recall a way the moduli spaces of
holomorphic curves can be oriented (see chapter 7 of
\cite{Bourgeois_thesis}) in case the asymptotics are non-degenerate. We
have a map $\pi$ from the moduli space of holomorphic curves $\mathcal M$
to the space of Fredholm operators $\mathcal O(\gamma^+,\gamma^-)$ with
behavior near the punctures corresponding to the closed Reeb orbits
$\gamma^+$ and $\gamma^-$ defined by sending a holomorphic map to its
linearized Cauchy-Riemann operator. Let $\mathcal L$ denote the determinant
bundle over $\mathcal O(\gamma^+,\gamma^-)$. Then $\pi^*\mathcal L$ is
naturally isomorphic to the top exterior power of $T\mathcal M$. In case
the asymptotics for $\mathcal O(\gamma^+,\gamma^-)$ are fixed, $\mathcal
O(\gamma^+,\gamma^-)$ is contractible and hence $\mathcal L$ is
trivial. In this case there are no obstructions to orient $\mathcal M$.

On the other hand, if the asymptotics are allowed to vary as is the case
for the Morse-Bott formalism, then we may get a non-contractible
space of Fredholm operators. This can happen, because the space of
operators fibers over submanifolds of the form $M_T$. If that space $M_T$
is not simply connected, it may contain a loop of operators such that the
determinant bundle over that loop is not trivial. Such a loop is called a {\bf
  disorienting loop}. We should remark here that if the projection of this
loop to the orbit space $S_T$ is contractible, the loop in $M_T$ is
homotopic to a bad Reeb orbit.

These phenomena can in general be present, because the linearized
Cauchy-Riemann operator is usually only real linear and not complex linear. In
favorable cases, the linearized Cauchy-Riemann operator is complex linear
as well and an orientation on the determinant line bundle can be obtained
directly from the induced complex structure on the kernel and cokernel
of the linearized Cauchy-Riemann operator. This removes the need to see
whether there are disorienting loops or bad orbits, because they cannot occur in that case.

\subsubsection{Differential for cylindrical contact homology in the Morse-Bott case}
With the above remarks in mind, the differential can be defined. The differential of the chain complex for cylindrical contact homology is given by
\begin{equation}
\label{MorseBott_differential}
dp=\partial p+\sum_q n_{q}q,
\end{equation}
where $p\in S_T$, $\partial$ is the Morse-Witten differential of the Morse
function $f_T$ on $S_{T'}$, $q\in S_{T'}$ and $n_q$ is the algebraic number of elements in
the fibered product
$$
(W^u(p)\times_{S_T} \mathcal M^{f_T}(S_T;S_{T'})\times_{S_{T'}}W^s(q))/\R
$$
if this product is zero-dimensional and $0$ otherwise ($q\in S_{T'}$). In this
product $W^s$ and $W^u$ denote the stable and unstable manifolds of a
critical point of $f_T$ on an orbit space, respectively. 

The precise definition of $n_q$ is more involved. However, we shall argue that in our case of a Brieskorn manifold and almost complex structure there are no zero-dimensional moduli spaces of holomorphic curves. As a result, the second term in the differential vanishes for us.

\begin{theorem}[Bourgeois]
\label{MorseBott_theorem}
Let $\alpha$ be a contact form of Morse-Bott type for a contact structure
$\xi$ on $M$ that satisfies $c_1(\xi)=0$. Assume that all Reeb orbits are closed. Assume that, for all $T\in
\sigma(\alpha)$, $M_T$ and $S_T$ are orientable, $\pi_1(S_T)$ has no
disorienting loop, and all Reeb orbits in $S_T$ are good. Assume that the
almost complex structure $J$ is invariant under the Reeb flow on all
submanifolds $M_T$. Assume that the cylindrical homology is well defined:
the Morse-Bott chain complex has no contractible orbits of index -1,0 or
+1. Assume furthermore that all orbit spaces $S_T$ of contractible periodic
orbits have index positivity or index negativity. Then the homology $H_*(C^{\bar a}_*)$ of the
Morse-Bott chain complex $(C^{\bar a}_*,d)$ of $(M,\alpha)$ is isomorphic to the cylindrical
homology $HF^{\bar a}_*(M,\xi)$ of $(M,\xi)$, where both homologies have coefficients in the
Novikov ring of $H_2(M,\Z)/\mathcal R$.
\end{theorem}

There are other (improved) versions of this theorem, but this one
suffices for our needs. In addition we will take the ring $\mathcal R$ in
the above theorem to be $H_2(M,\Z)$, or in other words we will use $\Q$
coefficients for the Morse-Bott chain complex.

\section{Algorithm for the computation of the cylindrical homology of
  some Brieskorn manifolds}

Consider the Brieskorn manifold $M=\Sigma(a_0,\ldots,a_n)\subset \C^{n+1}$ with
contact form induced by $\alpha=\frac{i}{8}\sum a_j(z_j \mathrm d \bar
z_j-\bar z_j \mathrm d z_j)$ and assume $n\geq 3$, which means that $M$ is at least
5-dimensional. The Reeb flow of the contact form $\alpha$ is
given by $\phi_t(z)=(e^{4it/a_0}z_0,\ldots,e^{4it/a_n}z_n)$. As we remarked in
the previous section, the first Chern class of the contact
structure induced by $\alpha$ is zero and all Reeb orbits are
closed. The description of the algorithm will use the notation introduced
in Section \ref{prelims}. The proof is in Section
\ref{proof_algorithm}. The algorithm works as follows.
\begin{itemize}
 \item[(1)] Compute the homology of $M$ using the algorithm of Randell
   \cite{Randell} which we described in Section
   \ref{Brieskorn_homology}. This information can be used to determine more
   precisely what manifold $M$ is (if the dimension of $M$ is five, this
   step provides enough information to use the classification of Barden,
   see \cite{Barden}). This step also involves some numerics that are used in
   Step $(4)$.

 \item[(2)] Identify all orbit types that can occur for the Reeb flow. This
   is done as follows. For
   all subsets $I_s\subset I=\{0,\ldots,n \}$, $s>1$, the minimal positive time $T$ such that
   $\frac{2T}{\pi}$ is divisible by all elements of $a_i$, $i\in I_s$, is
   $\frac{\pi}{2}\lcm_{i\in I_s}(a_i)$. The same minimal
   time $T$ can occur for several sets $I_s$. Let $J_T$ denote the largest
   such set. We
   obtain a collection of sets $J_{T_1},\ldots,J_{T_k}$ for different
   $T_1,\ldots,T_k$. The corresponding submanifolds $M_{T_i}:=K(J_{T_i})$ indicate
   submanifolds that are invariant under time $T_i$ of the Reeb flow.
 \item[(3)] Compute the Maslov indices of an orbit in $M_{T_i}$ with time
   $NT_i$ ($N\in\N$). In order to ensure that we do not consider orbits of
   another orbit space, we have to choose $N$ such that for $i\neq j$
   multiples $NT_i$ are not divisible by $T_j$ whenever $J_{T_i}\subset
   J_{T_j}$. Note that the Maslov index will only depend on the orbit type
   and not on the orbit itself. We may use the following formula
\begin{equation}
\label{eq_Maslov_index_Brieskorn}
\mu(S_{NT_i})=2\sum_{j \in J_{T_i}} \frac{2NT_i}{\pi a_j}+
2\sum_{j \in I-J_{T_i}} \lfloor \frac{2NT_i}{\pi a_j} \rfloor+\#(I-J_{T_i})
-4N\frac{T_i}{\pi}.
\end{equation}
The algorithm fails if one does not obtain index positivity or negativity
for the orbit spaces at this point. The conditions for these are given by 
\begin{equation*}
%%\label{index_pos_condition}
\sum_{j=0}^n\frac{1}{a_j}>1 \text{ for index positivity, and   }
\sum_{j=0}^n\frac{1}{a_j}<1 \text{ for index negativity.}
\end{equation*}
See also Remark \ref{rem_indexpos} for a discussion if these conditions are not met.

\item[(4)] The dimension of the orbit space $S_{T_i}=M_{T_i}/S^1$ is given by
  $2\# J_{T_i}-4$. Compute the rational homology of the orbit spaces
  $S_{T_i}$ in the following way. First of all, compute the rank of $\tilde H_{\# J_{T_i}-2}(M_{T_i})$, given by
$$
\kappa=\rk \tilde H_{\# J_{T_i}-2}(M_{T_i})=\sum_{I_s \subset J_{T_i}} (-1)^{\#
 J_{T_i} -s}\frac{\prod_{j'\in
    I_s}a_{j'}}{\underset{j\in I_s}{\lcm}~ a_j }. 
$$
By our remark in Section \ref{S1_remark} and Formula \ref{eq:rathom_orbit} we obtain
$$
\hspace{1 cm}
H_q(S_{T_i},\Q) \cong \left\{ \begin{array}{l} \Q, ~q \textrm{ even, } 0
    \leq q \leq \dim S_{T_i} \\ 0 \textrm{, otherwise} \end{array} \right\} \oplus
\left\{ \begin{array}{l} \Q^{\kappa},~q=\frac{1}{2} \dim S_{T_i} \\ 0 \textrm{ otherwise} \end{array}\right\}.
$$

\item[(5)] The cylindrical contact homology with $\Q$-coefficients of $M$
  with induced contact structure is a $\Q$-vector space, where the number
  of generators in each degree can be determined as follows. For each $T_i$
  we get $\rk H_j(S_{T_i},\Q)$ generators in degree
  $\mu(S_{NT_i})+n-3+j-\frac{1}{2}\dim S_{T_i} $ for
  $j=0,\ldots,\dim S_{T_i}$ and $N\in \N$ such that for $j\neq i$
  the multiples $NT_i$ are not divisible by $T_j$ whenever $J_{T_i}\subset
  J_{T_j}$. Here we use the Maslov-index that has been computed in Step 3.
\end{itemize}
For the cylindrical contact homology to be well defined and an invariant of
the contact structure there should be no
generators in degree $-1$, $0$ or $1$. To check this, one needs to define Morse
functions $f_T$ on the orbit spaces $S_T$ following Bourgeois's construction
described in Section \ref{MorseBott_homology}. The critical points of these
Morse functions form the Morse-Bott chain complex with grading given by
Formula \ref{eq:MorseBott_degree}. If there are no critical points with degree $-1,0$ or
$1$ then the algorithm yields the cylindrical contact homology for
contractible Reeb orbits. Note that these computations can depend on the
choice of Morse functions.
\begin{remark}
\label{periodic_CCH}
We should emphasize at this point that cylindrical contact homology of
Brieskorn manifolds is a periodic repetition of certain $\Q$ vector spaces
with a degree shift. This can be seen as follows. Let us consider the orbit space $S_{T_i}$ and multiple
coverings $S_{NT_i}$ where $N$ is chosen such that for $i\neq j$ the time
$NT_i$ is not divisible by $T_j$ whenever $J_{T_i}\subset J_{T_j}$. If we add 
$$
s_i:=\frac{\pi \underset{j'\in I}{\lcm} a_{j'}}{2 T_i}
$$
to $N$, we remain in an orbit space of the same type as
$\frac{\pi \lcm_{j' \in I} a_{j'}}{2}$ is divisible by all $T_j$. We see
that the Maslov index changes as follows,
$$
\mu(S_{(N+s_i)T_i})=\mu(S_{NT_i})+2~ \underset{j\in I}{\lcm} a_j ~\left(\sum_{j'=0}^n\frac{1}{a_{j'}}-1\right).
$$ 
This shift of the Maslov index is independent of the orbit space $S_{T_i}$
and hence the terms in the contact homology corresponding to the orbit
space $S_{T_i}$ are repeated with period at most 
$$
2~\underset{j\in I}{\lcm}~ a_j~\left( \sum_{j=0}^n\frac{1}{a_j}-1 \right).
$$
The periodicity of contact homology allows us to stop the algorithm after a
finite number of steps.
\end{remark}
\begin{remark}
Note that the requirement that $n\geq 3$, is not strictly necessary. For $n=2$ though, we are looking at 3-dimensional
Brieskorn manifolds, which are in general not simply connected. We can of
course deal with these cases in an easy way by considering only
contractible or homologically trivial Reeb orbits, but I have not worked out the details.
\end{remark}
\begin{remark}
If we consider Brieskorn manifolds with large exponents, we can ensure that
we have index negativity. In addition, large exponents ensure that the
grading is strictly less than $-1$, guaranteeing that cylindrical contact
homology is well defined. Indeed if
$$
\min_{i\in I}a_i\geq \frac{5n}{2},
$$
then the algorithm will always give the cylindrical contact homology of the
Brieskorn manifold $\Sigma(a_0,\ldots,a_n)$. This estimate is rather rough
and can be obtained from the formula for the Maslov index from Step (3).
\end{remark}
\begin{remark}
\label{rem_indexpos}
One might wonder what happens if index positivity/negativity fails. This is the case if
$$
\sum_{j=0}^n \frac{1}{a_j}=1.
$$
Note that in these cases there will always be degree $0$ orbits. Indeed, putting the full period $T=\pi/2\lcm_{i\in I}a_i$ into Formula \ref{eq_Maslov_index_Brieskorn} yields Maslov index $0$. By using the expression for the degree in part (5) of the algorithm we see that there will always be degree $0$ orbits. This means that, for our choice of contact form, cylindrical contact homology is not well defined as an invariant. 
\end{remark}

\subsection{Examples}

\subsubsection{Cylindrical homology of some contact structures on
  $S^2\times S^3$}
\label{contact_sigma(2l,2,2,2)}
In this section we consider the family of Brieskorn manifolds of the form
$M=\Sigma(2l,2,2,2)$ for $l>1$. Using our algorithm, it turns out that
these manifolds are diffeomorphic to $S^2\times S^3$ and that their
cylindrical contact homologies are all isomorphic. This is a bit
exceptional, since typically we get very different homologies for different
exponents. In the Section \ref{special_contact_sphere} we find some new exotic contact structures
on spheres which illustrates the latter point. Let us now turn our
attention to the application of the algorithm to $M=\Sigma(2l,2,2,2)$. The
numbering is as in the algorithm.
\begin{itemize}
\item[(1)] We find $\rk H_2(M)=1$. Computation of the
homology torsion by Randell's algorithm shows that there is none. By the
classification of simply connected five manifolds \cite{Barden} we see that
$M\cong S^2\times S^3$, as the second Stiefel-Whitney class of $M$
is zero.
\item[(2)] The time $T_1=\frac{\pi}{2}2$ is the minimal time for the sets
  which do not include $0$. Hence we see $J_{T_1}=\{ 1,2,3\}$.\\
  The time $T_2=\frac{\pi}{2}2l$ appears as minimal time for the set
  I=$\{ 0,1,2,3 \}$. So $J_{T_2}=\{ 0,1,2,3\}$.\\

\item[(3)] We get for $NT_1\frac{2}{\pi}$ not divisible by $l$
$$
\mu(S_{NT_1})=2N+1+ 2 \lfloor \frac{N}{l} \rfloor.
$$
The principal orbits have Maslov index
$$
\mu(S_{NT_2})=2lN+2N.
$$
\item[(4)] We find $\dim S_{T_1}=2$ with $H_0(S_{T_1},\Q)=\Q$ and
  $H_2(S_{T_1},\Q)=\Q$. All other homology groups are 0. The orbit space $S_{T_2}$ has dimension 4, with the homology
  $H_0(S_{T_2},\Q)=\Q$, $H_2(S_{T_2},\Q)=\Q^2$ and
  $H_4(S_{T_3},\Q)=\Q$. The other homology groups are zero. 
\item[(5)] The period we defined in Remark \ref{periodic_CCH} is equal to
  $2+2l$. This allows us to compute fewer terms. For the first period, the
  contributions from $S_{NT_1}$ lie in degree
$$
2N+2 \lfloor \frac{N}{l} \rfloor +k 
$$
for $N=1,\ldots,l-1$ (since we are considering a single period) and
$k=0,2$. The contribution due to $S_{NT_2}$ are in degree
$$
2l+k,
$$
for $k=0,2,4$. For the first period, we get one generator in degree
$2$, two generators in degree $4,6,\ldots,2l+2$ and one generator in degree
$2l+4$. The cylindrical contact homology has hence one generator in degree
$2$ and two generators in degree $2k$ for $k>1$. We note that the
cylindrical homology is well defined as there are no generators in degree
$-1,0$ or $1$ (lowest degree is higher than $1$).
\end{itemize}
We can also apply the algorithm to $\Sigma(2,2,2,2)$, where there is just a
single orbit space. This yields the same contact homology, see the
discussion in \ref{Lerman_discussion}.

\subsubsection{Some contact structures with index negativity}
\label{example_index_negativity}
Let us consider Brieskorn manifolds with large exponents such that we have
index negativity and that the degree is strictly less than $-1$. In case all exponents are equal, contact homology is particularly easy to
compute. For simplicity, we consider examples coming from Brieskorn
manifolds of the form $M=\Sigma(k,k,k,k)$ with $k\geq 7$. Since $M$ has only a single orbit space,
the computations are simple.

Step (1) shows that $\Sigma(k,k,k,k)$ is diffeomorphic to a connected sum
of $k^3-4k^2+6k-3$ copies of $S^2\times S^3$. The other steps yield the
following. The minimal return time is $T=\frac{\pi}{2}k$. Hence we get the Maslov
index
$$
\mu(S_{NT})=2\cdot 4\cdot N -2 \cdot N \cdot k=2N(4-k).
$$
Step (4) of the algorithm shows that $H_0(S_{T};\Q)\cong \Q$,
$H_2(S_{T};\Q)\cong \Q^d$, where $d=\kappa+1=\frac{(k-1)^4-1}{k}+2$. The
last homology group is $H_4(S_{T};\Q)\cong \Q$. This gives a single period
of contact homology. Taking the Maslov index into account we get the
following for the contact homology of $M$.

We have one generator in degree $2N(4-k)-2$ for $N=1,2,\ldots$. In degree
$2N(4-k)$ we have $\frac{(k-1)^4-1}{k}+2$ generators and in degree
$2N(4-k)+2$ we have again one generator ($N=1,2,\ldots$).

\subsection{Proof of the algorithm}
\label{proof_algorithm}
Consider the Brieskorn manifold given by $M=\Sigma(a_0,\ldots,a_n)\subset \C^{n+1}$ with
contact form induced by $\alpha=\frac{i}{8}\sum a_j(z_j \mathrm d \bar
z_j-\bar z_j \mathrm d z_j)$. The differential $\mathrm d \alpha=\frac{i}{4}\sum
a_j \mathrm d z_j \w \mathrm d \bar z_j$ is a symplectic form on
$\C^{n+1}$. Let us denote this symplectic form by $\omega=\mathrm d\alpha$. Let $\xi$ be the
contact structure given by $\ker \alpha|_M$. Note that the contact
form is of Morse-Bott type (Definition \ref{MorseBott_definition}). This is seen as follows. Discreteness of the
action spectrum is guaranteed by Step (2) of the algorithm. The sets
$M_T$ in Definition \ref{MorseBott_definition} are Brieskorn manifolds with
their standard contact form. In particular they are closed submanifolds of
$M$. Note that this verifies the rank condition on $\mathrm d \alpha$ as well. The
last condition for $\alpha$ being of Morse-Bott type can be checked by
observing that the differential of the Reeb flow $\phi$ is diagonal, see
Formula \ref{eq:Reebflow} which allows us to verify the rank condition on
$T\phi -id$.

We verify that Step (3) gives the correct Maslov indices. For each time $T_i$
that we found in Step (2) consider the $N$-fold covering of $M_{T_i}$ with
$N$ such that $NT_i$ is not divisible by $T_j$'s of a larger orbit space.
Now let $p\in M_{T_i}$ and consider its orbit under the Reeb flow
$\gamma(t):=\phi_t(p)$ for $t\in [0,NT_i]$. To compute the indices we
use Formula~\ref{eq:maslov_reeb}. We denote the obvious extension of the
Reeb flow to $\C^{n+1}$ also by $\phi_t$.
 
The symplectic action of the extended Reeb flow $\phi_t$ on $\C^{n+1}$ is given by the differential of
$T\phi_t$, the diagonal matrix $\diag(e^{4it/a_0},\ldots,e^{4it/a_n})$ ($t\in[0,NT_i]$). As in Section \ref{maslov_section} we denote the path of
symplectic matrices induced by the extended flow by $\Phi_{\C^{n+1}}$. We
can now use the additivity of the Maslov index and Formula~\ref{eq:maslov_unitair} to get the index of this path.  This gives us
$$
\mu(\Phi_{\C^{n+1}})=2\sum_{a_j \in J_{T_i}} \frac{2NT_i}{\pi a_j}+
2\sum_{a_j \in I-J_{T_i}} \lfloor \frac{2NT_i}{\pi a_j} \rfloor+\#(I-J_{T_i}).
$$

With respect to the symplectic basis of $\xi^\omega$ found in
Section \ref{prelims} the action induced on $\xi^\omega$ by the extended
Reeb flow is given by Formula \ref{eq:maslov_complement}. As before we use
the notation from Section \ref{maslov_section}, i.e. $\Phi_{\xi^\omega}$
denotes the path of symplectic matrices induced by the action of the
extended Reeb flow. Formula \ref{eq:maslov_unitair} gives us
$$
\mu(\Phi_{\xi^\omega})=4N\frac{T_i}{\pi}.
$$
Taking the difference of these Maslov indices yields the desired result
from Step (3). The conditions for index positivity and negativity can be
found by observing that $\lfloor t \rfloor \geq t-1$ and $\lfloor t \rfloor \leq t$.

We show that the determinant bundle of the linearized Cauchy-Riemann
operator can be oriented directly, i.e.~we shall show that the
linearized Cauchy-Riemann operator is asymptotic to a complex linear operator. First of all, we note
that near the punctures this operator can be given the form (see for
instance \cite{Bourgeois_thesis}, chapter 3 and 5)
$$
\frac{\partial}{\partial s}+J_0 \frac{\partial}{\partial t} +S(s,t),
$$
where $S$ are symmetric matrices and $J_0$ is the standard complex
structure and $(s,t)$ are cylindrical coordinates near the puncture, $t$
for the $S^1$ direction. We have 
$$
\lim_{s\to\infty} S(s,t)=\tilde S(t), 
$$
with
$$
\tilde S(t)=-J_0 \frac{d \psi(t)}{dt} \psi^{-1}(t),
$$
where $\psi$ are the symplectic matrices obtained from the linearized Reeb
flow in the symplectization $M\times \R$. Note that if these matrices
$\psi(t)$ are unitary with respect to the above trivialization, the matrix
$\tilde S$ will commute with $J_0$. We will verify that we can choose a
trivialization such that the matrices $\psi$ are unitary.

Note that the linearized extended Reeb flow on $\C^{n+1}$ is represented by
a path of unitary matrices, see Section \ref{maslov_section}. Keep in mind that the metric is given by
$\omega(\ldots,J\ldots)$, where $J$ is the standard complex structure on
$\C^{n+1}$. For the above purposes we need a trivialization that comes from
the symplectization of $M$, so this one is not suitable. To that end,
recall that a contact structure $\xi$ is symplectically stably trivial if
and only if $\xi \oplus \C$ is symplectically trivial (see also Remark \ref{stably_trivial}). This means that we
can split off a complex line bundle from our trivialization on $\C^{n+1}$.
Define $V:=\mathrm{span}_{\R}(\tilde X_1,\tilde Y_1)$ using the
notation from Section \ref{prelims} and let $W$ be the
orthogonal complement of $V$ in $\C^{n+1}$. Because the linearized extended
Reeb
flow maps $V$ to $V$ (see Section \ref{maslov_section}), we can see that
the linearized extended Reeb flow induces
a map from $W$ to
$W$ which is unitary with respect to the induced metric. Note also that
$\xi$ is a subbundle of $W$ and $W$ contains the Reeb line bundle. Hence
$W$ can be identified with the tangent bundle of the symplectization of $M$
restricted to $M\times \{\text{point}\}$. Let $\psi(t)$ be the path of
unitary matrices given by the linearized extended Reeb flow on $W$ with
respect to an orthonormal basis of $W$. 

If we define $\tilde S(t)$ as above in terms of the unitary $\psi(t)$, then
it follows that the Cauchy-Riemann operator
$$
\frac{\partial}{\partial s}+J \frac{\partial}{\partial t} +\tilde S(t)
$$
is complex linear. By the above discussion, this operator is asymptotic
to the given linearized Cauchy-Riemann operator on the symplectization. As
mentioned in Section \ref{MorseBott_homology}, this gives us an orientation
of the determinant bundle. In particular, no bad orbits or disorienting
loops can occur.

Index positivity/negativity is verified in Step (3) of the algorithm. If we have that
there are no generators of degree $-1,0$ or $1$, then the homology of the
Morse-Bott complex is isomorphic to cylindrical contact homology according
to Theorem \ref{MorseBott_theorem}. Consider the Morse-Bott complex with
generators the critical points $p$ of the chosen Morse functions with
grading given by $\mu(S_T)-\frac{1}{2}\dim S_T+\mathrm{index}~p+n-3$.
The differential of the Morse-Bott complex is given by Formula
\ref{MorseBott_differential}.  The differential acting on $p\in S_T$ is
given by
$$
d p =\partial p+\sum n_q q.
$$
The first term is the Morse-Witten differential for the critical points of
the Morse functions $f_T$. The second term counts the number of elements in the zero-dimensional part of the fibered product
$$
(W^u(p)\times_S \mathcal M^{f_T}(S;S')\times_{S'}W^s(q))/\R.
$$
Now note that there is an $S^1$-action on the symplectization of $M$
induced by the Reeb flow. A holomorphic curve asymptotic to closed Reeb
orbits comes therefore in at least an $S^1$-family (by letting the Reeb flow
act) except in the case that the curve is a vertical cylinder. This means the
above fibered product is at least 1 dimensional, so it will not contribute
to the differential. This argument is also used in Section~9.3 of \cite{Bourgeois_thesis}.

This implies that the only non-zero contribution to the differential comes from $\partial p$, which
means that the cylindrical contact homology is isomorphic to the
Morse-Witten homology of the orbit spaces $S_T$ with degree shifts of
$\mu(S_T)-\frac{1}{2}\dim S_T+n-3$.
As the Morse-Witten homology of the orbit spaces is equal to the rational
homology of the orbit spaces which is computed in Step (4), we find that the contact homology is given by our algorithm.

\section{Exotic contact structures}
Our aim in this section is to describe a certain class of contact manifolds
that admits infinitely many non-isomorphic contact structures. 

Given two contact manifolds $(M_1,\xi_1)$ and $(M_2,\xi_2)$, we can build a
new contact manifold by forming their connected sum, \cite{Meckert} and
\cite{Weinstein}. If we think of a connected sum as first removing a disk from
both $M_1$ and $M_2$ and then gluing them via a ``connecting'' tube, then
the contact structure on $M_1\# M_2$ can be be made to coincide with the
contact structure on $M_1$ with a disk removed and $M_2$ with a disk removed.

In order to say something about the contact homology of the connected sum,
we find generic contact forms $\alpha_1,~\alpha_2$ for the contact
structures $\xi_1$ and $\xi_2$, i.e.~contact forms whose closed
Reeb orbits are non-degenerate.  

First we need another, but similar notion of index positivity, which we
will refer to as Ustilovsky index positivity. Suppose that a contact
structure $\xi$ is symplectically stably trivial and let $F$ be a
corresponding trivialization. We may then compute the Maslov index with
respect to the trivialization $F$. The index does depend on the
trivialization and we indicate this by making the Maslov index visibly
dependent on the trivialization $F$ by writing $\mu(\ldots;F)$. We use $\epsilon$ to indicate the trivial bundle.
\begin{definition} Let $(M,\alpha)$ be a contact manifold. Assume that
  $\pi_1(M)=0$ and that the bundle $(\xi,d\alpha)$ is symplectically stably trivial. Let
  $F$ be a symplectic trivialization of $\xi\oplus \epsilon^2$. The manifold
  $(M,\alpha)$ is called {\bf Ustilovsky index-positive} if there exist
  constants $c>0$ and $d$ such that for any Reeb orbit $\gamma$
$$
\mu(\gamma ;F)\geq c \mathcal A(\gamma)+d
$$
holds.
\end{definition}

\begin{remark}
\label{stably_trivial}
Ustilovsky has shown in his thesis that a contact structure $\xi$ is
symplectically stably
trivial if and only if $\xi\oplus \epsilon^2$ is trivial and that the
notion of Ustilovsky index positivity does not depend on the choice of
trivialization. Note that the Chern class of a symplectically stably
trivial contact structure is trivial.
\end{remark}
\begin{example}
\label{Ustilovsky_index_positivity}
Brieskorn manifolds are symplectically stably trivial (for instance, see
Section \ref{prelims}). Hence we may consider Ustilovsky index
positivity of a Brieskorn manifold $\Sigma(a_0,\ldots,a_n)$. Note that a Brieskorn manifold has
Ustilovsky index positivity if all orbits types have index
positivity. Namely, our
computation of the Maslov indices for closed orbits used the trivialization
coming from $\C^{n+1}$. If use that trivialization to verify Ustilovsky
index positivity, we see that the formula from Step (3) of the algorithm can be
modified as follows
$$
\mu(\gamma(t)|_{t\in[0,T]})=
\sum_{j=0}^n \mu(e^{4it/a_j}|_{t\in[0,T]})-\mu(e^{4it}|_{t\in[0,T]})
$$
where the Reeb orbit $\gamma$ is given by
$$
\gamma(t)=(e^{4it/a_0}z_0,\ldots,e^{4it/a_n}z_n),
$$
for $(z_0,\ldots,z_n)\in \Sigma(a_0,\ldots,a_n)\subset \C^{n+1}$. Then apply Formula \ref{eq:maslov_unitair}.
\end{example}

We have the following theorem from Ustilovsky (\cite{Ustilovsky_thesis}, Theorem 5.2.1):
\begin{theorem}[Ustilovsky] \label{connected_sum}
Let $(M_1,\alpha_1)$, $(M_2,\alpha_2)$ be two simply connected contact
manifolds of dimension $2n-1$ that have Ustilovsky index positivity. Assume all periodic
orbits of $R_{\alpha_1}$ and $R_{\alpha_2}$ are non-degenerate. Then
for any integer $N$ there exists a contact form $\alpha$ on $M=M_1 \# M_2$
so that
\begin{itemize}
\item[(1)] ($M,\alpha$) is Ustilovsky index-positive.
\item[(2)] All periodic orbits of $R_\alpha$ are non-degenerate.
\item[(3)] If $c_j^1$, $c_j^2$ and $c_j$ denote the numbers of periodic
  Reeb orbits of degree $j$ in $M_1$, $M_2$ and $M$, respectively, then
  for $j\leq N$, we have $c_j=c_j^1+c_j^2+\beta_j$, where $\beta_j=1$ for
  $j=2n-3,2n-1,\ldots$ and $\beta_j=0$ otherwise.
\end{itemize} 
\end{theorem}
The $\beta_j$'s in this theorem are the degrees of the periodic Reeb orbits
in the connecting tube. We will take for $M_1$ any contact manifold
satisfying the conditions of the above theorem of Ustilovsky and for $M_2$
we will take a special contact sphere. This sphere will have the property
that its contact homology contains generators with degree lower than the
lowest degree of a generator in the connecting tube. After taking the
connected sum with $M_1$ and using Ustilovsky's theorem, the resulting
contact manifold $M_1\# M_2$ will be diffeomorphic to $M_1$, but its
cylindrical contact homology will have more generators in low degrees.

\subsection{Construction of a special contact
  sphere}\label{special_contact_sphere}

Let us consider the Brieskorn manifold $M=\Sigma(p_1,\ldots,p_{n-1},2,2)$ where $p_1,\ldots,p_{n-1}$ are odd primes which
will be specified later (they need to be chosen large enough). Notice
that we immediately see that this manifold is homeomorphic to a sphere by
Theorem 14.5 in \cite{Hirzebruch}. This fact can also be seen by computing
the homology of $M$ using Randell's algorithm. Then use
the fact that Brieskorn manifolds of dimension at least 5 are always simply connected and conclude that $M$ is homeomorphic to a sphere with the
generalized Poincaré conjecture as proved by Smale. We apply our algorithm to compute the first terms of the
contact homology.

\subsubsection{Contact homology of
  $\Sigma(p_1,\ldots,p_{n-1},2,2)$}\label{contact_homology_sphere}

We have the following orbit types and Maslov indices. 
\begin{itemize}
\item $I_2=\{ 2,2 \}$. The minimal return time for orbits of this type is $T=\pi$. The
  Maslov index of the corresponding orbits is given by
$$
\mu(S_{NT})=2\sum_{i=1}^{n-1}\lfloor \frac{2N}{p_i} \rfloor+n-1 \geq n-1.
$$ 
If the $p_i$'s are odd primes, the first term will vanish for at least
$N=1$. Now we turn our attention to the homology of the orbit space. It has
dimension 0 and Formula \ref{eq:rathom_orbit} shows that $H_0(S_T,\Q)=\Q^2$. This
shows that the contact homology of $M$ has at least 2 generators in
degree $2n-4$. Since the first term is always even, multiple covers of this
orbit will have either the same degree or a higher even degree. Note that
this orbit type will not give any generators in degree $2n-3$. 
\item Sets of the form $\{ p_{i_1},\ldots,p_{i_k} \}$ with $k$ at least 2. The minimal return time is
  now $T=p_{i_1}\cdot \ldots \cdot p_{i_k}\frac{\pi}{2}$. The Maslov
indices of the orbit spaces are given by
\begin{eqnarray*}
\hspace{1 cm}\mu(S_{NT})&=& 2\sum_j p_{i_1}\cdot \ldots \cdot {\hat p}_{i_j}\cdot \ldots
\cdot p_{i_k}N+4 \lfloor
\frac{Np_{i_1}\cdot \ldots \cdot p_{i_k}}{2} \rfloor \\
\hspace{1 cm}
&+&2\sum_{l\neq i_j} \lfloor
\frac{Np_{i_1}\cdot \ldots \cdot p_{i_k}}{p_l} \rfloor +n-1-k+2
-2Np_{i_1}\cdot \ldots \cdot p_{i_k}.
\end{eqnarray*}
Note that $4 \lfloor \frac{Np_{i_1}\cdot \ldots \cdot p_{i_k}}{2} \rfloor\geq
2Np_{i_1}\cdot \ldots \cdot p_{i_k}-2$ and that $\dim S_T=2k-4$. Using this estimate,
we find that the degree of the associated generators can be estimated from below as
$$
\mathrm{degree} \geq 2n-2-2k+2\sum_j p_{i_1}\cdot \ldots \cdot \hat p_{i_j}
\cdot \ldots \cdot p_{i_k}N.
$$
Since the sum contains at least two terms, we can make the Maslov
index arbitrarily large by choosing big primes. In particular, the degree of
these orbit types can be assumed to be larger than $2n-3$.
\item Sets of the form $\{ p_{i_1},\ldots ,p_{i_k},2,2 \}$. The minimal return
  time is $T=p_{i_1}\cdot \ldots \cdot p_{i_k}\pi$. The associated Maslov indices
  are
$$
\hspace{1.25 cm}
\mu(S_{NT})=4\sum_j
p_{i_1}\ldots {\hat p}_{i_j}\cdot \ldots \cdot p_{i_k}N+2\sum_{l\neq i_j} \lfloor
\frac{2Np_{i_1}\cdot \ldots \cdot p_{i_k}}{p_l} \rfloor +n-1-k.
$$
For $k>1$ the first term will contain at least one $p_i$-term. This means
that the degree (note that the dimension of the orbit space is now $2k$)  will become as large as we like by
choosing the $p$'s accordingly. For $k=1$ the first term is $4N$ and we see
that the degree is at least $4N+n-3+n-3\geq 2n-2$.
\end{itemize}
Summarizing these estimates, we see that by choosing suitable primes, we
may assume that the contact homology contains at least two generators in
degree $2n-4$ and no generators in degree $2n-3$. Note that we have
(Ustilovsky) index positivity since two exponents are 2, see the condition
from Step (3) of the algorithm.

In order to be able to apply Theorem \ref{connected_sum}, we need to have a generic contact form on $M$ that has at least two
closed Reeb orbits in degree $2n-4$ and no generators in degree $2n-3$. We
will do this, along with a more general statement in the following
interlude.

\subsubsection{ Generic contact forms }

In this section, we want to associate a generic contact form, i.e.~a
contact form whose closed Reeb orbits are non-degenerate, to the
Morse-Bott contact forms used in our algorithm.

Let $M$ be a contact manifold of dimension $2n-1$ and let $\alpha$ be a
contact form on $M$ that satisfies the Morse-Bott condition. We use
Bourgeois's construction of Morse functions on the orbit spaces described in Section \ref{MorseBott_homology}. We use his ideas to perturb
$\alpha$ into a generic contact form such that we still have some
information on the indices of the closed Reeb orbits. The following
observation by Ustilovsky \cite{Ustilovsky_thesis} will play a key role.

\begin{lemma}\label{Ustilovsky_pos}
If $(M,\alpha)$ is Ustilovsky index-positive, then a $C^{\infty}$-small
perturbation $\alpha '$ of $\alpha$ the manifold $(M,\alpha ')$ is also
Ustilovsky index-positive. Moreover, if $\mu(\gamma ;F)\geq c \mathcal
A(\gamma)+d$ for orbits $\gamma$ of $R_\alpha$ then, for $\alpha'$ close
enough to $\alpha$, $\mu(\gamma ';F')\geq c' \mathcal A(\gamma ')+d'$ for
orbits $\gamma '$ of $R_{\alpha '}$, where $c'=c/2$ and $d'=-|d|-2n$.
\end{lemma}
\begin{remark}\label{Perturb_remark}
Note that for a small perturbation, we have a one-to-one correspondence between
non-degenerate Reeb orbits of the contact form and of the perturbed
contact form up to some period. The Conley-Zehnder indices of the
corresponding orbits are the same.
\end{remark}
\begin{lemma} \label{Perturb_lemma} Let $\alpha$ be the standard contact form on the Brieskorn
  manifold $M=\Sigma(a_0,\ldots ,a_n)$, $n\geq 3$. Assume that the exponents
  are such that we have Ustilovsky index positivity (cf.~Example \ref{Ustilovsky_index_positivity}). Then for all $N\in \Z$ there exists
  a generic contact form $\alpha '$, such that all generators of the chain complex of $\alpha '$ coincide with the generators of the Morse-Bott chain complex
  of $\alpha$ up to degree $N$.
\end{lemma}
\begin{proof} Let us denote the constants from the Ustilovsky index positivity
by $c>0$ and $d$ such that $\mu(\gamma)\geq c\mathcal A(\gamma)+d$ for a
part $\gamma$ of a Reeb orbit. Choose $T\geq
\mathrm{max}\{ N,\frac{4}{c}(N+|d|+4n) \}$. By Lemma \ref{Bourgeois_perturb} we find a
perturbation $\alpha ''$ of $\alpha$ such that its periodic Reeb orbits
up to action $T$ are non-degenerate and correspond to critical points of
the chosen Morse functions on the orbit spaces. Since we have Ustilovsky
index positivity for $\alpha$, we will have the same for $\alpha ''$ by Lemma \ref{Ustilovsky_pos}, where the constants are now $c/2$
and $-|d|-2n$. 

We perturb $\alpha ''$ further to make all Reeb orbits
non-degenerate and call the perturbation $\alpha '$. Once again $\alpha '$
is Ustilovsky index positive with constants $c/4$ and $-|d|-4n$. By Remark
\ref{Perturb_remark} this perturbation will not change the Conley-Zehnder
indices of orbits with period up to $T$ (and in particular up to $N$). The perturbation will
introduce new periodic Reeb orbits. However, the newly created ones can be
made to have period larger than $T$ (by making a small enough
perturbation). As a result of Lemma \ref{Ustilovsky_pos}, their Conley-Zehnder
indices will be larger than $T c/4 -|d|-4n\geq N$. This proves our lemma.
\end{proof}

\subsubsection{A generic contact form for $\Sigma(p_1,\ldots ,p_{n-1},2,2)$}
We continue with our construction of a contact sphere from Section
\ref{special_contact_sphere}. As before, we write our Brieskorn sphere as
$M=\Sigma(p_1,\ldots ,p_{n-1},2,2)$. We get the Morse-Bott chain complex for
contact homology by following Bourgeois's construction of Morse functions
for the orbit spaces from Section \ref{MorseBott_homology}. Note that the
computation of the contact homology of $M$ in Section
\ref{contact_homology_sphere} shows that the Morse-Bott complex will have
at least two generators in degree $2n-4$ and no generators in degree $2n-3$.
This holds true because the lowest dimensional orbit spaces have dimension
0, so the number and degree of generators associated to those orbit spaces
do not depend on the choice of Morse functions. Note that these data do
depend on the choice of Morse functions for the other orbit spaces.

Now we apply Lemma \ref{Perturb_lemma} to obtain a generic contact form
where the the number and degree of generators of the chain complex coincide
with those of the Morse-Bott complex up to degree $2n-2$. This gives $M$ a
generic contact form and allows us to use Theorem \ref{connected_sum}.

\subsection{Constructing new contact structures}

\begin{theorem}\label{contact_structures}
Let $(M,\xi)$ be a simply-connected contact manifold. Assume furthermore that $M$
admits a nice contact form (a contact form without any closed Reeb orbits
of degree $-1,0$ or $1$) that has Ustilovsky index positivity (in
particular $c_1(\xi)=0$). Then $M$ admits infinitely many non-isomorphic contact structures.
\end{theorem}

\begin{proof} Let $N'=\Sigma (p_1,\ldots ,p_{n-1},2,2)$. The above discussion
shows that $N'$ admits a generic contact form that is Ustilovsky index
positive with at least two generators in degree $2n-4$, no generators in
degree $2n-3$ and no generators in lower degrees. As remarked before, we
know that $N'$ is homeomorphic to a sphere (Theorem 14.5 in
\cite{Hirzebruch}). Since the differentiable structures on a sphere form a
finite group, we can find an $r$ such that $N:=\underbrace{N'\#\ldots \#N'}_{r}$ is
diffeomorphic to the standard sphere. By Theorem \ref{connected_sum}, $N$
admits a generic contact form whose cylindrical contact homology has at least $2r$ generators in degree $2n-4$,
precisely $r-1$ generators in degree $2n-3$ and no generators in lower degrees.

Now apply Theorem \ref{connected_sum} to $M$ and $N$. The connected sum $M\# N$ will be diffeomorphic
to $M$, since $N$ is diffeomorphic to $S^{2n-1}$. The theorem shows that the
connected sum still admits a nice contact form (because no generators are added
in degrees $-1,0$ or $1$). The number of generators of the chain complex of
the cylindrical contact homology is increased by at least $2r$ in degree $2n-4$ and by
$r$ in degree $2n-3$. Because of Ustilovsky index positivity there are only a finite
number of generators in each degree. Since the number of generators in
degree $2n-5$ is unchanged by taking the connected sum with $N$, taking
repeated connected sums with $N$ will ensure that we get a contact
structure on $M$ with a contact homology different from the original one,
namely with more generators in degree
$2n-4$. By taking more connected sums with $N$ we get infinitely many
contact structures on $M$, all with different cylindrical contact homology distinguished by the rank of the homology in degree $2n-4$.
\end{proof}

\begin{remark}
Note that if two contact manifolds satisfy the conditions of Theorem
\ref{contact_structures}, their connected sum satisfies those conditions as
well by Theorem \ref{connected_sum}. In addition, any index positive
Brieskorn manifold whose cylindrical contact homology can be computed using our algorithm
satisfies the conditions of Theorem \ref{contact_structures}. This can be
seen by applying Lemma \ref{Perturb_lemma}.

As an example, the cylindrical contact homology of
$\Sigma(2,\ldots ,2)=ST^*S^n$ can be computed using our algorithm (the
computation is similar to the one in Section \ref{example_index_negativity}), which
shows that $\#_k ST^*S^n$ carries infinitely many contact structures for all
$k\in \N$ and $n>2$.
\end{remark}

\begin{question}
Consider the orbifold $X_a=\Sigma(a_0,\ldots,a_n)/S^1$. We can compute the
orbifold cohomology of $X_a$ following Chen and Ruan, see
\cite{Chen_Ruan}. This is done by taking the direct sum of the homologies
of the singular subspaces with appropriate degree shifts. This is similar
to what we do in our algorithm, where our ``degree shifts'' are expressed
in the Maslov indices.

Orbifold cohomology is not quite the same as contact homology, though.
First of all, contact homology is periodic, whereas orbifold cohomology has
only non-trivial groups in a finite number of degrees. Secondly, the
orbifold cohomology can have rational grading and contact homology has
integer grading in our cases.

A single period of contact homology could however be related to orbifold
cohomology. Indeed, in case $X_a$ is a manifold and not just an orbifold (for
instance when all exponents of the Brieskorn manifold are equal), contact
homology can be expressed as a periodic repetition of the homology of $X_a$. In general, the degree shifts of orbifold cohomology do not match
the ones coming from the Maslov index. We can, of course, transform the
degree shifts to correct this. Can this be done in a meaningful way? This
might give some insight of how to compute contact homology for more general
$S^1$-bundles over symplectic orbifolds.
\end{question}

\begin{question}
\label{Lerman_discussion}
In \cite{Lerman} Lerman gives a description of a number of contact structures
on $ST^*S^3 \cong S^2\times S^3$. Some are distinguished by their first Chern class, but
he is also able to produce many contact structures with the same Chern
class. One can wonder whether these contact structures are isomorphic. The
contact structures we found in \ref{contact_sigma(2l,2,2,2)} are related to this example. Indeed,
Lerman produces infinitely many contact structures in each even first Chern
class except when $c_1$ is zero. This is complemented by the Brieskorn
manifolds $\Sigma(2k,2,2,2)$ which all have trivial Chern class. Note
however that the cylindrical contact homology is not able to distinguish
between those contact structures. This leaves the question whether one is
able to distinguish these possibly non-isomorphic contact structures in
each even Chern class.
\end{question}

\begin{question}
Can one show the result of Theorem \ref{contact_structures} under
milder conditions? In particular, to what extend is index positivity necessary?
\end{question}

\bibliographystyle{plain}

\end{document}